\documentclass[leqno,12pt]{amsart}
\usepackage{amssymb} 
\theoremstyle{plain} 
\newtheorem{theorem}{\indent\sc Theorem}[section] 
\newtheorem{lemma}[theorem]{\indent\sc Lemma}

\newtheorem{proposition}[theorem]{\indent\sc Proposition}

\theoremstyle{definition} 
\newtheorem{definition}[theorem]{\indent\sc Definition}

\begin{document}

\title[partial differential systems in two variables]{Five-parameter family of partial differential systems in two variables \\}
\author{Yusuke Sasano }

\renewcommand{\thefootnote}{\fnsymbol{footnote}}
\footnote[0]{2000\textit{ Mathematics Subjet Classification}.
34M55; 34M45; 58F05; 32S65.}

\maketitle

\begin{abstract}
We find a five-parameter family of partial differential systems in two variables with two polynomial Hamiltonians. We give its symmetry and holomorphy conditions. These symmetries, holomorphy conditions and invariant divisors are new.
\end{abstract}

\section{Introduction}
In this paper, we present a 5-parameter family of partial differential systems in two variables explicitly given by
\begin{align}\label{1}
\begin{split}
dq_1&=\frac{\partial H_1}{\partial p_1}dt+\frac{\partial H_2}{\partial p_1}ds, \quad dp_1=-\frac{\partial H_1}{\partial q_1}dt-\frac{\partial H_2}{\partial q_1}ds,\\
dq_2&=\frac{\partial H_1}{\partial p_2}dt+\frac{\partial H_2}{\partial p_2}ds, \quad dp_2=-\frac{\partial H_1}{\partial q_2}dt-\frac{\partial H_2}{\partial q_2}ds
\end{split}
\end{align}
with the polynomial Hamiltonians:
\begin{align}\label{2}
\begin{split}
H_1 &=H_{VI}(q_1,p_1,t;\alpha_1,\alpha_2,\alpha_3,\alpha_4)\\
&+\alpha_2 p_2\left\{\frac{-(t-1)sq_1+t(s-1)q_2+(t-s)q_1q_2}{t(t-1)(t-s)}+\frac{(q_1-t)q_2(q_2-1)}{t(t-1)(t-\eta)}\right\}\\
&+\alpha_5 p_1\left\{\frac{(t-s)q_1(q_1-1)+t(t-1)(q_1-q_2)}{t(t-1)(t-s)}+\frac{(q_1-t)((t-1)q_1+(q_1-t)q_2)}{t(t-1)(t-\eta)}\right\}\\
&- p_1p_2\left\{\frac{(t-1)(sq_1^2+tq_2^2)-(t-s)q_2(q_1^2+t)-2t(s-1)q_1q_2}{t(t-1)(t-s)}-\frac{(q_1-t)^2q_2(q_2-1)}{t(t-1)(t-\eta)}\right\}\\
&+\frac{\alpha_2\alpha_5(2tq_1-q_1-tq_2+q_1q_2-\eta q_1)}{t(t-1)(t-\eta)},\\
H_2&=\pi(H_1),
\end{split}
\end{align}
where the transformation $\pi$ is explicitly given by
\begin{align}
\begin{split}
\pi:&(q_1,p_1,q_2,p_2,t,s;\alpha_0,\alpha_1,\alpha_2,\alpha_3,\alpha_4,\alpha_5)\\
&\rightarrow(q_2,p_2,q_1,p_1,s,t;\alpha_0,\alpha_1,\alpha_5,\alpha_3,\alpha_4,\alpha_2).
\end{split}
\end{align}
Here $q_1,p_1,q_2$ and $p_2$ denote unknown complex variables, and $\alpha_0,\alpha_1,\ldots,\alpha_5$ are complex parameters satisfying the relation:
\begin{equation}\label{3}
\alpha_0+\alpha_1+2\alpha_2+\alpha_3+\alpha_4+2\alpha_5=1.
\end{equation}
This parameter's relation can be obtained by holomorphy conditions in Theorem \ref{th;holoG}.

The symbol $H_{VI}(q,p,t;\beta_1,\beta_2,\beta_3,\beta_4)$ denotes the Hamiltonian of the second-order Painlev\'e VI equations (see \cite{Sasa5}) given by
\begin{align}\label{6}
\begin{split}
&t(t-1)(t-\eta)H_{VI}(q,p,t;\beta_1,\beta_2,\beta_3,\beta_4)\\
&=q(q-1)(q-\eta)(q-t)p^2\\
&+\{\beta_1(t-\eta)q(q-1)+2\beta_2q(q-1)(q-\eta)\\
&+\beta_3(t-1)q(q-\eta)+\beta_4 t(q-1)(q-\eta)\}p\\
&+\beta_2\{(\beta_1+\beta_2)(t-\eta)+\beta_2(q-1)\\
&+\beta_3(t-1)+t \beta_4\}q \quad (\beta_0+\beta_1+2\beta_2+\beta_3+\beta_4=1, \quad \eta \in {\Bbb C}-\{0,1\}).
\end{split}
\end{align}

We give its symmetry and holomorphy conditions. These symmetries, holomorphy conditions and invariant divisors are new.

After we review the notion of accessible singularity and local index, we make its holomorphy conditions by resolving the accessible singularities.

\section{Symmetry and holomorphy conditions}
In this section, we give its symmetry and holomorphy conditions. These symmetries, holomorphy conditions and invariant divisors are new.
\begin{theorem}\label{th;holoG}
Let us consider a polynomial Hamiltonian system with Hamiltonians $H_i \in {\Bbb C}(t,s)[q_1,p_1,q_2,p_2] \ (i=1,2)$. We assume that

$(A1)$ $deg(H_i)=6$ with respect to $q_1,p_1,q_2,p_2$.

$(A2)$ This system becomes again a polynomial Hamiltonian system in each coordinate $r_i, \ i=0,1,\dots,5${\rm:\rm}
\begin{align}
\begin{split}
&r_0:x_0=-p_1((q_1-t)p_1+(q_2-s)p_2-\alpha_0), \ y_0=\frac{1}{p_1}, \ z_0=(q_2-s)p_1, \ w_0=\frac{p_2}{p_1}, \\
&r_1:x_1=-p_1((q_1-\eta)p_1+(q_2-\eta)p_2-\alpha_1), \ y_1=\frac{1}{p_1}, \ z_1=(q_2-\eta)p_1, \ w_1=\frac{p_2}{p_1}, \\
&r_2:x_2=\frac{1}{q_1}, \ y_2=-q_1(q_1p_1+\alpha_2), \ z_2=q_2, \ w_2=p_2, \\
&r_3:x_3=-p_1((q_1-1)p_1+(q_2-1)p_2-\alpha_3), \ y_3=\frac{1}{p_1}, \ z_3=(q_2-1)p_1, \ w_3=\frac{p_2}{p_1}, \\
&r_4:x_4=-p_1(q_1p_1+q_2p_2-\alpha_4), \ y_4=\frac{1}{p_1}, \ z_4=q_2p_1, \ w_4=\frac{p_2}{p_1}, \\
&r_5:x_5=q_1, \ y_5=p_1, \  z_5=\frac{1}{q_2}, \ w_5=-(q_2p_2+\alpha_5)q_2.
\end{split}
\end{align}
Then such a system coincides with the system \eqref{1} with two polynomial Hamiltonians \eqref{2}.
\end{theorem}

\begin{proposition} 
In each coordinate $r_i, \ i=0,1,\dots,5$, the  Hamiltonians $H_{j1}$ and $H_{j2}$ on $U_j \times B$ are expressed as a polynomial in $x_j,y_j,z_j,w_j$ and a rational function in $t$ and $s$, and satisfy the following conditions{\rm: \rm}
\begin{align}\label{symplectic}
\begin{split}
&dq_1 \wedge dp_1 +dz \wedge dp_2 - dH_1 \wedge dt- dH_2 \wedge ds\\
&=dx_j \wedge dy_j +dz_j \wedge dw_j - dH_{j1} \wedge dt- dH_{j2} \wedge ds \quad (j=1,2,\dots,5),\\
&dq_1 \wedge dp_1 +dq_2 \wedge dp_2 - d(H_1-p_1) \wedge dt- d(H_2-p_2) \wedge ds\\
&=dx_0 \wedge dy_0 +dz_0 \wedge dw_0 - dH_{01} \wedge dt- dH_{02} \wedge ds.
\end{split}
\end{align}
\end{proposition}

\begin{figure}
\unitlength 0.1in
\begin{picture}(35.87,16.56)(17.90,-24.16)
%
\special{pn 20}%
\special{ar 2062 1574 267 244  0.0000000 6.2831853}%
%
\special{pn 20}%
\special{ar 3110 1602 267 245  0.0000000 6.2831853}%
%
\special{pn 20}%
\special{ar 4110 1608 267 244  0.0000000 6.2831853}%
%
\special{pn 20}%
\special{ar 5110 1615 267 243  0.0000000 6.2831853}%
%
\special{pn 20}%
\special{ar 3582 906 192 146  0.0000000 6.2831853}%
%
\special{pn 20}%
\special{ar 3598 2270 192 146  0.0000000 6.2831853}%
%
\special{pn 8}%
\special{pa 3390 930}%
\special{pa 2206 1352}%
\special{fp}%
%
\special{pn 8}%
\special{pa 3470 1026}%
\special{pa 3174 1363}%
\special{fp}%
%
\special{pn 8}%
\special{pa 3702 1026}%
\special{pa 4022 1369}%
\special{fp}%
%
\special{pn 8}%
\special{pa 3774 934}%
\special{pa 4966 1391}%
\special{fp}%
%
\special{pn 8}%
\special{pa 2198 1796}%
\special{pa 3406 2218}%
\special{fp}%
%
\special{pn 8}%
\special{pa 3206 1836}%
\special{pa 3446 2149}%
\special{fp}%
%
\special{pn 8}%
\special{pa 4942 1813}%
\special{pa 3782 2218}%
\special{fp}%
%
\special{pn 8}%
\special{pa 3966 1813}%
\special{pa 3718 2144}%
\special{fp}%
%
\special{pn 8}%
\special{pa 1790 1574}%
\special{pa 2318 1574}%
\special{fp}%
%
\special{pn 8}%
\special{pa 2854 1586}%
\special{pa 3358 1586}%
\special{fp}%
%
\special{pn 8}%
\special{pa 3846 1597}%
\special{pa 4358 1597}%
\special{fp}%
%
\special{pn 8}%
\special{pa 4854 1608}%
\special{pa 5366 1608}%
\special{fp}%
\put(35.0200,-9.4700){\makebox(0,0)[lb]{$p_1$}}%
\put(35.2600,-23.4400){\makebox(0,0)[lb]{$p_2$}}%
\put(18.5400,-15.2300){\makebox(0,0)[lb]{$q_1-\eta$}}%
\put(18.6200,-17.2200){\makebox(0,0)[lb]{$q_2-\eta$}}%
\put(29.7400,-15.3400){\makebox(0,0)[lb]{$q_1$}}%
\put(29.7400,-17.2200){\makebox(0,0)[lb]{$q_2$}}%
\put(39.4200,-17.1600){\makebox(0,0)[lb]{$q_2-s$}}%
\put(49.2600,-15.5100){\makebox(0,0)[lb]{$q_1-1$}}%
\put(49.3400,-17.1600){\makebox(0,0)[lb]{$q_2-1$}}%
\put(39.4200,-15.4500){\makebox(0,0)[lb]{$q_1-t$}}%
\end{picture}%
\label{UraDynkin1}
\caption{The symbol in each circle denotes the invariant cycle for the system.}
\end{figure}
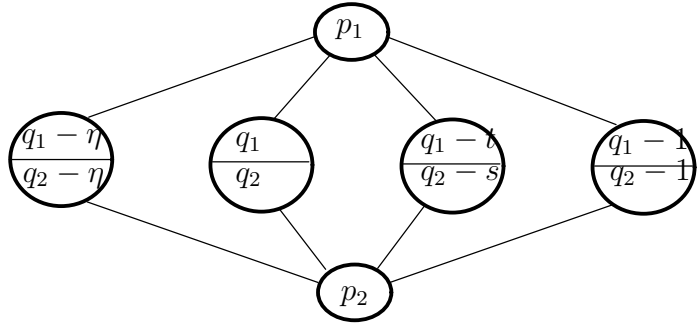

\begin{center}\label{inv}
\begin{tabular}{|c|c|c|} \hline
codimension & invariant cycles & parameter's relation \\ \hline
1 & $f_2:=p_1$ & $\alpha_2=0$  \\ \hline
1 & $f_5:=p_2$ & $\alpha_5=0$  \\ \hline
2 & $f_0^{(1)}:=q_1-t, \ f_0^{(2)}:=q_2-s$ & $\alpha_0=0$  \\ \hline
2 & $f_1^{(1)}:=q_1-\eta, \ f_1^{(2)}:=q_2-\eta$ & $\alpha_1=0$   \\ \hline
2 & $f_3^{(1)}:=q_1-1, \ f_3^{(2)}:=q_2-1$ & $\alpha_3=0$   \\ \hline
2 & $f_4^{(1)}:=q_1, \ f_4^{(2)}:=q_2$ & $\alpha_4=0$  \\ \hline
\end{tabular}
\end{center}
We note that when $\alpha_2=0$, we see that the system \eqref{1} admits a particular solution $f_2=0$, and when $\alpha_0=0$, we see that the system \eqref{1} admits a particular solution $f_0^{(1)}=f_0^{(2)}=0$.

\section{B{\"a}cklund transformations}
\begin{theorem}\label{th:in3}
The system \eqref{1} admits the following transformations as its B{\"a}ckl-\\
und transformations\rm{:\rm} with the notation $(*)=(q_1,p_1,q_2,p_2,\eta,t,s;\alpha_0,\alpha_1,\dots,\alpha_5),$
\begin{align}
\begin{split}
        s_1: (*) \rightarrow &\left(q_1+\frac{\alpha_2}{p_1},p_1,q_2,p_2,\eta,t,s;\alpha_0+\alpha_2,\alpha_1+\alpha_2,-\alpha_2,\alpha_3+\alpha_2,\alpha_4+\alpha_2,\alpha_5 \right), \\
        s_2: (*) \rightarrow &\left(q_1,p_1,q_2+\frac{\alpha_5}{p_2},p_2,\eta,t,s;\alpha_0+\alpha_5,\alpha_1+\alpha_5,\alpha_2,\alpha_3+\alpha_5,\alpha_4+\alpha_5,-\alpha_5 \right),\\
        \pi_1: (*) \rightarrow &(\frac{\eta-q_1}{\eta-1},-(\eta-1)p_1,\frac{\eta-q_2}{\eta-1},-(\eta-1)p_2,\frac{\eta}{\eta-1},\frac{\eta-t}{\eta-1},\frac{\eta-s}{\eta-1};\\
        &\alpha_0,\alpha_4,\alpha_2,\alpha_3,\alpha_1,\alpha_5),\\
        \pi_2: (*) \rightarrow &(\frac{q_1(t-\eta)}{t-q_1+tq_1-\eta t},-\frac{(t-q_1+tq_1-\eta t)\{(t-q_1+tq_1-\eta t)p_1+\alpha_2(t-1)\}}{t(t-\eta)(\eta-1)},\\
        &\frac{q_2(s-\eta)}{s-q_2+sq_2-\eta s},-\frac{(s-q_2+sq_2-\eta s)\{(s-q_2+sq_2-\eta s)p_2+\alpha_5(s-1)\}}{s(s-\eta)(\eta-1)},\\
        &\eta,\frac{\eta-t}{1-2t+\eta t},\frac{\eta-s}{1-2s+\eta s};\alpha_3,\alpha_1,\alpha_2,\alpha_0,\alpha_4,\alpha_5),\\
        \pi_3: (*) \rightarrow &(\frac{(t-1)q_1}{t-q_1-\eta t+\eta tq_1},\frac{(t-q_1+\eta t(q_1-1))\{(q_1-t)p_1+\alpha_2-\eta t((q_1-1)p_1+\alpha_2)\}}{t(t-1)(\eta-1)},\\
        &\frac{(s-1)q_2}{s-q_2-\eta s+\eta sq_2},\frac{(s-q_2+\eta s(q_2-1))\{(q_2-s)p_2+\alpha_5-\eta s((q_2-1)p_2+\alpha_5)\}}{s(s-1)(\eta-1)},\\
        &\frac{1}{\eta},\frac{\eta(t-1)}{t-\eta-\eta t+\eta^2 t},\frac{\eta(s-1)}{s-\eta-\eta s+\eta^2 s};\alpha_1,\alpha_0,\alpha_2,\alpha_3,\alpha_4,\alpha_5),\\
        \pi_4: (*) \rightarrow &(1-q_1,-p_1,1-q_2,-p_2,1-\eta,1-t,1-s;\alpha_0,\alpha_1,\alpha_2,\alpha_4,\alpha_3,\alpha_5),\\
        \pi_5: (*) \rightarrow &(q_2,p_2,q_1,p_1,\eta,s,t;\alpha_0,\alpha_1,\alpha_5,\alpha_3,\alpha_4,\alpha_2).
        \end{split}
        \end{align}
\end{theorem}
The B{\"a}cklund transformations $s_1,s_2$ are determined by the invariant divisors \eqref{inv}.

\section{Accessible singularity and local index}
Let us review the notion of {\it accessible singularity}. Let $B$ be a connected open domain in $\Bbb C$ and $\pi : {\mathcal W} \longrightarrow B$ a smooth proper holomorphic map. We assume that ${\mathcal H} \subset {\mathcal W}$ is a normal crossing divisor which is flat over $B$. Let us consider a rational vector field $\tilde v$ on $\mathcal W$ satisfying the condition
\begin{equation*}
\tilde v \in H^0({\mathcal W},\Theta_{\mathcal W}(-\log{\mathcal H})({\mathcal H})).
\end{equation*}
Fixing $t_0 \in B$ and $P \in {\mathcal W}_{t_0}$, we can take a local coordinate system $(x_1,\ldots ,x_n)$ of ${\mathcal W}_{t_0}$ centered at $P$ such that ${\mathcal H}_{\rm smooth \rm}$ can be defined by the local equation $x_1=0$.
Since $\tilde v \in H^0({\mathcal W},\Theta_{\mathcal W}(-\log{\mathcal H})({\mathcal H}))$, we can write down the vector field $\tilde v$ near $P=(0,\ldots ,0,t_0)$ as follows:
\begin{equation*}
\tilde v= \frac{\partial}{\partial t}+g_1 
\frac{\partial}{\partial x_1}+\frac{g_2}{x_1} 
\frac{\partial}{\partial x_2}+\cdots+\frac{g_n}{x_1} 
\frac{\partial}{\partial x_n}.
\end{equation*}
This vector field defines the following system of differential equations
\begin{equation}\label{39}
\frac{dx_1}{dt}=g_1(x_1,\ldots,x_n,t),\ \frac{dx_2}{dt}=\frac{g_2(x_1,\ldots,x_n,t)}{x_1},\cdots, \frac{dx_n}{dt}=\frac{g_n(x_1,\ldots,x_n,t)}{x_1}.
\end{equation}
Here $g_i(x_1,\ldots,x_n,t), \ i=1,2,\dots ,n,$ are holomorphic functions defined near $P=(0,\dots ,0,t_0).$

\begin{definition}\label{Def1}
With the above notation, assume that the rational vector field $\tilde v$ on $\mathcal W$ satisfies the condition
$$
(A) \quad \tilde v \in H^0({\mathcal W},\Theta_{\mathcal W}(-\log{\mathcal H})({\mathcal H})).
$$
We say that $\tilde v$ has an {\it accessible singularity} at $P=(0,\dots ,0,t_0)$ if
$$
x_1=0 \ {\rm and \rm} \ g_i(0,\ldots,0,t_0)=0 \ {\rm for \rm} \ {\rm every \rm} \ i, \ 2 \leq i \leq n.
$$
\end{definition}

If $P \in {\mathcal H}_{{\rm smooth \rm}}$ is not an accessible singularity, all solutions of the ordinary differential equation passing through $P$ are vertical solutions, that is, the solutions are contained in the fiber ${\mathcal W}_{t_0}$ over $t=t_0$. If $P \in {\mathcal H}_{\rm smooth \rm}$ is an accessible singularity, there may be a solution of \eqref{39} which passes through $P$ and goes into the interior ${\mathcal W}-{\mathcal H}$ of ${\mathcal W}$.

Here we review the notion of {\it local index}. Let $v$ be an algebraic vector field with an accessible singular point $\overrightarrow{p}=(0,\ldots,0)$ and $(x_1,\ldots,x_n)$ be a coordinate system in a neighborhood centered at $\overrightarrow{p}$. Assume that the system associated with $v$ near $\overrightarrow{p}$ can be written as
\begin{align}\label{b}
\begin{split}
\frac{d}{dt}\begin{pmatrix}
             x_1 \\
             x_2 \\
             \vdots\\
             x_{n-1} \\
             x_n
             \end{pmatrix}=\frac{1}{x_1}\left\{\begin{bmatrix}
             a_{11} & 0 & 0 & \hdots & 0 \\
             a_{21} & a_{22} & 0 &  \hdots & 0 \\
             \vdots & \vdots & \ddots & 0 & 0 \\
             a_{(n-1)1} & a_{(n-1)2} & \hdots & a_{(n-1)(n-1)} & 0 \\
             a_{n1} & a_{n2} & \hdots & a_{n(n-1)} & a_{nn}
             \end{bmatrix}\begin{pmatrix}
             x_1 \\
             x_2 \\
             \vdots\\
             x_{n-1} \\
             x_n
             \end{pmatrix}+\begin{pmatrix}
             x_1h_1(x_1,\ldots,x_n,t) \\
             h_2(x_1,\ldots,x_n,t) \\
             \vdots\\
             h_{n-1}(x_1,\ldots,x_n,t) \\
             h_n(x_1,\ldots,x_n,t)
             \end{pmatrix}\right\},\\
              (h_i \in {\Bbb C}(t)[x_1,\ldots,x_n], \ a_{ij} \in {\Bbb C}(t))
             \end{split}
             \end{align}
where $h_1$ is a polynomial which vanishes at $\overrightarrow{p}$ and $h_i$, $i=2,3,\ldots,n$ are polynomials of order at least 2 in $x_1,x_2,\ldots,x_n$, We call ordered set of the eigenvalues $(a_{11},a_{22},\cdots,a_{nn})$ {\it local index} at $\overrightarrow{p}$.

We are interested in the case with local index
\begin{equation}\label{integer}
(1,a_{22}/a_{11},\ldots,a_{nn}/a_{11}) \in {\Bbb Z}^{n}.
\end{equation}
These properties suggest the possibilities that $a_1$ is the residue of the formal Laurent series:
\begin{equation}
y_1(t)=\frac{a_{11}}{(t-t_0)}+b_1+b_2(t-t_0)+\cdots+b_n(t-t_0)^{n-1}+\cdots \quad (b_i \in {\Bbb C}),
\end{equation}
and the ratio $(1,a_{22}/a_{11},\ldots,a_{nn}/a_{11})$ is resonance data of the formal Laurent series of each $y_i(t) \ (i=2,\ldots,n)$, where $(y_1,\ldots,y_n)$ is original coordinate system satisfying $(x_1,\ldots,x_n)=(f_1(y_1,\ldots,y_n),\ldots,f_n(y_1,\ldots,y_n)), \ f_i(y_1,\ldots,y_n) \in {\Bbb C}(t)(y_1,\ldots,y_n)$.

If each component of $(1,a_{22}/a_{11},\ldots,a_{nn}/a_{11})$ has the same sign, we may resolve the accessible singularity by blowing-up finitely many times. However, when different signs appear, we may need to both blow up and blow down.

The $\alpha$-test,
\begin{equation}\label{poiuy}
t=t_0+\alpha T, \quad x_i=\alpha X_i, \quad \alpha \rightarrow 0,
\end{equation}
yields the following reduced system:
\begin{align}\label{ppppppp}
\begin{split}
\frac{d}{dT}\begin{pmatrix}
             X_1 \\
             X_2 \\
             \vdots\\
             X_{n-1} \\
             X_n
             \end{pmatrix}=\frac{1}{X_1}\begin{bmatrix}
             a_{11}(t_0) & 0 & 0 & \hdots & 0 \\
             a_{21}(t_0) & a_{22}(t_0) & 0 &  \hdots & 0 \\
             \vdots & \vdots & \ddots & 0 & 0 \\
             a_{(n-1)1}(t_0) & a_{(n-1)2}(t_0) & \hdots & a_{(n-1)(n-1)}(t_0) & 0 \\
             a_{n1}(t_0) & a_{n2}(t_0) & \hdots & a_{n(n-1)}(t_0) & a_{nn}(t_0)
             \end{bmatrix}\begin{pmatrix}
             X_1 \\
             X_2 \\
             \vdots\\
             X_{n-1} \\
             X_n
             \end{pmatrix},
             \end{split}
             \end{align}
where $a_{ij}(t_0) \in {\Bbb C}$. Fixing $t=t_0$, this system is the system of the first order ordinary differential equation with constant coefficient. Let us solve this system. At first, we solve the first equation:
\begin{equation}
X_1(T)=a_{11}(t_0)T+C_1 \quad (C_1 \in {\Bbb C}).
\end{equation}
Substituting this into the second equation in \eqref{ppppppp}, we can obtain the first order linear ordinary differential equation:
\begin{equation}
\frac{dX_2}{dT}=\frac{a_{22}(t_0) X_2}{a_{11}(t_0)T+C_1}+a_{21}(t_0).
\end{equation}
By variation of constant, in the case of $a_{11}(t_0) \not= a_{22}(t_0)$ we can solve explicitly:
\begin{equation}
X_2(T)=C_2(a_{11}(t_0)T+C_1)^{\frac{a_{22}(t_0)}{a_{11}(t_0)}}+\frac{a_{21}(t_0)(a_{11}(t_0)T+C_1)}{a_{11}(t_0)-a_{22}(t_0)} \quad (C_2 \in {\Bbb C}).
\end{equation}
This solution is a single-valued solution if and only if
$$
\frac{a_{22}(t_0)}{a_{11}(t_0)} \in {\Bbb Z}.
$$
In the case of $a_{11}(t_0)=a_{22}(t_0)$ we can solve explicitly:
\begin{equation}
X_2(T)=C_2(a_{11}(t_0)T+C_1)+\frac{a_{21}(t_0)(a_{11}(t_0)T+C_1){\rm Log}(a_{11}(t_0)T+C_1)}{a_{11}(t_0)} \quad (C_2 \in {\Bbb C}).
\end{equation}
This solution is a single-valued solution if and only if
$$
a_{21}(t_0)=0.
$$
Of course, $\frac{a_{22}(t_0)}{a_{11}(t_0)}=1 \in {\Bbb Z}$.
In the same way, we can obtain the solutions for each variables $(X_3,\ldots,X_n)$. The conditions $\frac{a_{jj}(t)}{a_{11}(t)} \in {\Bbb Z}, \ (j=2,3,\ldots,n)$ are necessary condition in order to have the Painlev\'e property.

\section{Construction of the holomorphy conditions}
In this section, we will give the holomorphy conditions $r_i \ (i=0,1,\ldots,5)$ by resolving some accessible singular loci of the system \eqref{1}.

In order to consider the singularity analysis for the system \eqref{1}, as a compactification of ${\Bbb C}^4$ which is the phase space of the system \eqref{1}, we take 4-dimensional complex manifold $\mathcal S$ given in the paper \cite{Sasa5}. This manifold can be considered as a generalization of the Hirzebruch surface.

We easily see that the rational vector field $\tilde v$ associated with the system \eqref{1} satisfies the condition:
$$
\tilde v \in H^0({\mathcal S},\Theta_{\mathcal S}(-\log{\mathcal H})({\mathcal H})).
$$

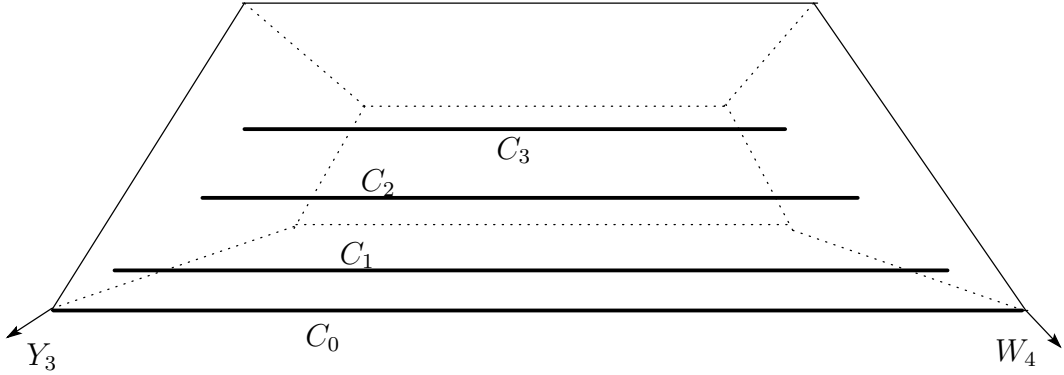
\begin{figure}
\unitlength 0.1in
\begin{picture}( 55.1000, 18.0000)( 15.7000,-28.0000)
%
\special{pn 8}%
\special{pa 2810 1000}%
\special{pa 1800 2610}%
\special{fp}%
\special{pa 1800 2610}%
\special{pa 6910 2610}%
\special{fp}%
%
\special{pn 8}%
\special{pa 2800 1000}%
\special{pa 5810 1000}%
\special{fp}%
%
\special{pn 8}%
\special{pa 5790 1000}%
\special{pa 6900 2610}%
\special{fp}%
%
\special{pn 8}%
\special{pa 2810 1000}%
\special{pa 3440 1540}%
\special{dt 0.045}%
\special{pa 3440 1540}%
\special{pa 3070 2180}%
\special{dt 0.045}%
%
\special{pn 8}%
\special{pa 1800 2600}%
\special{pa 3070 2170}%
\special{dt 0.045}%
%
\special{pn 8}%
\special{pa 3430 1540}%
\special{pa 5320 1540}%
\special{dt 0.045}%
\special{pa 5320 1540}%
\special{pa 5790 1000}%
\special{dt 0.045}%
%
\special{pn 8}%
\special{pa 5330 1530}%
\special{pa 5680 2190}%
\special{dt 0.045}%
\special{pa 5680 2190}%
\special{pa 6890 2610}%
\special{dt 0.045}%
%
\special{pn 8}%
\special{pa 3080 2160}%
\special{pa 5650 2160}%
\special{dt 0.045}%
%
\special{pn 8}%
\special{pa 1800 2600}%
\special{pa 1570 2750}%
\special{fp}%
\special{sh 1}%
\special{pa 1570 2750}%
\special{pa 1638 2730}%
\special{pa 1616 2722}%
\special{pa 1616 2698}%
\special{pa 1570 2750}%
\special{fp}%
%
\special{pn 8}%
\special{pa 6900 2610}%
\special{pa 7080 2800}%
\special{fp}%
\special{sh 1}%
\special{pa 7080 2800}%
\special{pa 7050 2738}%
\special{pa 7044 2762}%
\special{pa 7020 2766}%
\special{pa 7080 2800}%
\special{fp}%
\put(16.7000,-29.2000){\makebox(0,0)[lb]{$Y_3$}}%
\put(67.4000,-29.0000){\makebox(0,0)[lb]{$W_4$}}%
%
\special{pn 20}%
\special{pa 1810 2610}%
\special{pa 6880 2610}%
\special{fp}%
%
\special{pn 20}%
\special{pa 2130 2400}%
\special{pa 6490 2400}%
\special{fp}%
%
\special{pn 20}%
\special{pa 2590 2020}%
\special{pa 6020 2020}%
\special{fp}%
%
\special{pn 20}%
\special{pa 2810 1660}%
\special{pa 5640 1660}%
\special{fp}%
\put(31.3000,-28.2000){\makebox(0,0)[lb]{$C_0$}}%
\put(33.1000,-23.9000){\makebox(0,0)[lb]{$C_1$}}%
\put(34.2000,-20.1000){\makebox(0,0)[lb]{$C_2$}}%
\put(41.3000,-18.4000){\makebox(0,0)[lb]{$C_3$}}%
\end{picture}%
\label{Garnierfig1}
\caption{This figure denotes the boundary divisor ${\mathcal H}$ of ${\mathcal S}$. The bold lines $C_i \ i=0,1,2,3$ in ${\mathcal H}$ denote the accessible singular loci of the system \eqref{1}.}
\end{figure}

\begin{lemma}\label{lem1}
The rational vector field $\tilde v$ has the following accessible singular loci \rm{(see figure 2)}$:$
\begin{equation}
  \left\{
  \begin{aligned}
    C_0=&\{(X_3,Y_3,Z_3,W_3)|X_3=t,Z_3=s,Y_3=0\}\\
&\cup \{(X_4,Y_4,Z_4,W_4)|X_4=t,Z_4=s,W_4=0\} \cong {\Bbb P}^1,\\
    C_1=&\{(X_3,Y_3,Z_3,W_3)|X_3=\eta,Z_3=\eta,Y_3=0\}\\
&\cup \{(X_4,Y_4,Z_4,W_4)|X_4=\eta,Z_4=\eta,W_4=0\} \cong {\Bbb P}^1,\\
    C_2=&\{(X_3,Y_3,Z_3,W_3)|X_3=1,Z_3=1,Y_3=0\}\\
&\cup \{(X_4,Y_4,Z_4,W_4)|X_4=1,Z_4=1,W_4=0\} \cong {\Bbb P}^1,\\
    C_3=&\{(X_3,Y_3,Z_3,W_3)|X_3=Z_3=Y_3=0\}\\
&\cup \{(X_4,Y_4,Z_4,W_4)|X_4=Z_4=W_4=0\} \cong {\Bbb P}^1.
   \end{aligned}
  \right. 
\end{equation}
\end{lemma}
Here, the coordinate systems $(X_i,Y_i,Z_i,W_i) \ (i=3,4)$ (see \cite{Sasa5}) are explicitly given by
\begin{align}
\begin{split}
(X_3,Y_3,Z_3,W_3)=\left(q_1,\frac{1}{p_1},q_2,\frac{p_2}{p_1} \right),\\
(X_4,Y_4,Z_4,W_4)=\left(q_1,\frac{p_1}{p_2},q_2,\frac{1}{p_2} \right).
\end{split}
\end{align}

This lemma can be proven by a direct calculation.

Next, we calculate its local index at the point $P:=\{(X_3,Y_3,Z_3,W_3)|X_3=t,Z_3=s,Y_3=W_3=0\}$.

{\bf Step 0:} We make a change of variables.
$$
X_3^{(1)}=X_3-t, \quad Y_3^{(1)}=Y_3, \quad Z_3^{(1)}=Z_3-s, \quad W_3^{(1)}=W_3.
$$
Around the point $P$, we rewrite the system \eqref{1} as follows:
\begin{align*}
\frac{d}{dt}\begin{pmatrix}
             X_3^{(1)} \\
             Y_3^{(1)} \\
             Z_3^{(1)} \\
             W_3^{(1)}
             \end{pmatrix}&=\frac{1}{Y_3^{(1)}}\left\{\begin{pmatrix}
             2 & -\alpha_0 & 0 & 0  \\
             0 & 1 & 0 & 0 \\
             0 & 0 & 1 & 0 \\
             0 & \frac{\alpha_5}{t-s} & 0 & 0
             \end{pmatrix}\begin{pmatrix}
             X_3^{(1)} \\
             Y_3^{(1)} \\
             Z_3^{(1)} \\
             W_3^{(1)} 
             \end{pmatrix}+\cdots\right\}.
             \end{align*}
We see that this system has its local index $(2,1,1,0)$ at the point $P$.

For the remaining accessible singular loci, the local index is same.

\begin{proposition}\label{prop3}
If we resolve the accessible singular loci given in Lemma \ref{lem1} by blowing-ups, then we can obtain the canonical coordinate systems $r_i \ (i=0,1,3,4)$.
\end{proposition}

{\it Proof.} By the following steps, we can resolve the accessible singular locus $C_0$.

{\bf Step 1:} We blow up along the curve $C_0$.
$$
X_3^{(2)}=\frac{X_3^{(1)}}{Y_3^{(1)}}, \quad Y_3^{(2)}=Y_3^{(1)}, \quad Z_3^{(2)}=\frac{Z_3^{(1)}}{Y_3^{(1)}}, \quad W_3^{(2)}=W_3^{(1)}.
$$

{\bf Step 2:} We blow up along the surface $\{(X_3^{(2)},Y_3^{(2)},Z_3^{(2)},W_3^{(2)})|X_3^{(2)}=-Z_3^{(2)} W_3^{(2)}+\alpha_0\}$
$$
X_3^{(3)}=\frac{X_3^{(2)}+Z_3^{(2)} W_3^{(2)}-\alpha_0}{Y_3^{(2)}}, \quad Y_3^{(3)}=Y_3^{(2)}, \quad Z_3^{(3)}=Z_3^{(2)}, \quad W_3^{(3)}=W_3^{(2)}.
$$
By choosing a new coordinate system as
$$
(x_0,y_0,z_0,w_0)=(-X_3^{(3)},Y_3^{(3)},Z_3^{(3)},W_3^{(3)}),
$$
we can obtain the coordinate system $r_0$.

For the remaining accessible singular loci, the proof is similar. Thus, we have completed the proof of Proposition \ref{prop3}.

\begin{proposition}
After a series of explicit blowing-ups given in Proposition \ref{prop3}, we obtain the smooth projetive $4$-fold $\tilde{\mathcal S}$ and a birational morphism $\varphi:\tilde{\mathcal S} \rightarrow {\mathcal S}$. Its canonical divisor $K_{\tilde{\mathcal S}}$ is given by
\begin{equation}
K_{\tilde{\mathcal S}}=-3\tilde{\mathcal H}-\sum_{i=0}^3{\mathcal E}_i,
\end{equation}
where the symbol $\tilde{\mathcal H}$ denotes the proper transform of ${\mathcal H}$ by $\varphi$ and ${\mathcal E}_i$ denote the exceptional divisors obtained by Step $1$ \rm{(see Proof of Proposition \ref{prop3})}.
\end{proposition}

\end{document}